\documentclass{article}

\usepackage{epic}
\usepackage{eepic}
\usepackage{amssymb}
\usepackage{latexsym}
\usepackage{theorem}

\newcommand{\be}{\begin{equation}}
\newcommand{\ee}{\end{equation}} 
\newcommand{\bea}{\begin{eqnarray}}
\newcommand{\beann}{\begin{eqnarray*}}
\newcommand{\eea}{\end{eqnarray}}
\newcommand{\eeann}{\end{eqnarray*}}

\newcommand{\al}{\alpha}

\makeatletter
\@addtoreset{equation}{section}
\makeatother

\theoremstyle{break}

\newtheorem{p1}{Proposition}[section]
\newtheorem{l1}{Lemma}[section]

\newtheorem{thp1}[p1]{Theorem}
\newtheorem{thl1}[l1]{Theorem}

\newtheorem{pp1}[p1]{Proposition}

\newtheorem{cp1}[p1]{Corollary}

\begin{document}

\begin{center}
{\Large\bf 
Automorphic forms, fake monster algebras \\
and hyperbolic reflection groups}\\[1cm]
Nils R. Scheithauer\footnote{Supported by the Emmy-Noether program.},\\
Department of Mathematics,\\
University of California,
Berkeley, CA 94720\\
nrs@math.berkeley.edu\\
\end{center}
\vspace*{2cm}

\noindent
We construct 2 families of automorphic forms related to twisted fake monster algebras and calculate their Fourier expansions. This gives a new proof of their denominator identities and shows that they define automorphic forms of singular weight. We also obtain new infinite product identities which are the denominator identities of generalized Kac-Moody superalgebras. Finally we describe the reflection groups of the root lattices of these algebras.


\section{Introduction}

Borcherds, Gritsenko and Nikulin have shown that there are interesting relations between automorphic forms, generalized Kac-Moody algebras and hyperbolic reflection groups. They can be summarized roughly as follows. First the denominator identities of nice generalized Kac-Moody algebras often define automorphic forms. Second many reflection groups of Lorentzian lattices are associated to automorphic forms whose singularities are at the reflection hyperplanes of the reflection group. Finally the root lattices of nice generalized Kac-Moody algebras often have nice hyperbolic reflection groups. In all 3 cases it is not known what the precise necessary conditions are. For example it is not known when the denominator function of a generalized Kac-Moody algebra is an automorphic form. 

In this paper we give some new examples for the above relations. 
We construct 2 families of automorphic forms of singular weight. The automorphic forms are the denominator functions of generalized Kac-Moody algebras similar to the fake monster algebras. In the bosonic case the reflection groups of their root lattices are associated to automorphic forms. Then we show that the reflection groups of their root lattices are similar to those of $II_{1,25}$ and $II_{1,9}$.

We describe the sections of this paper in more detail. 

In the second section we recall some results about lattices and the singular theta correspondence. 

We use this correspondence in section $3$ to construct a family of automorphic forms related to twisted fake monster algebras. We give a new proof of their denominator identities and show that they define automorphic forms. 

In section $4$ we derive analogous results for twisted fake monster superalgebras. In additon we get new infinite product identities which are the denominator identities of generalized Kac-Moody superalgebras. 

In the last section we describe the reflection groups of the root lattices of the generalized Kac-Moody algebras in this paper. In the bosonic case they are related to automorphic forms. In the super case we use Vinberg's algorithm to determine their fundamental domains.

\section{Lattices and automorphic forms}

In this section we fix some notations and recall some results on lattices and the singular theta correspondence \cite{B3}.

\subsection{Lattices}

Let $M$ be a lattice with dual $M'$. We write $M(n)$ for the lattice obtained from $M$ by multiplying all norms with $n$. The level of $M$ is the smallest positive integer $n$ so that $n\lambda^2\in 2{\mathbb Z}$ for all $\lambda\in M'$. It follows $nM'\subset M$. 
We consider some examples. 

The lattice $II_{1,1}(n)$ consists of the elements $(m_1,m_2)\in {\mathbb Z}^2$ of norm $(m_1,m_2)^2
\linebreak
=-2nm_1m_2$. The lattice has determinant $n^2$ and the quotient of the dual $II_{1,1}(n)'$ by $II_{1,1}(n)$ is ${\mathbb Z}_n^2$. We write the elements of $II_{1,1}(n)'$ as $(m_1/n,m_2/n)$ with $m_i\in {\mathbb Z}$ so that $(m_1/n,m_2/n)^2=-2m_1m_2/n$. This shows that $II_{1,1}(n)$ has level $n$. Clearly $II_{1,1}(n)'$ is isomorphic to $II_{1,1}(1/n)$. 

Let $E_8^p$ be the sublattice of $E_8$ fixed by an automorphism of cycle shape $1^mp^m$ where $p$ is a prime such that $m=8/(p+1)$ is an integer. Then $E_8^p$ has level $p$ and determinant $p^m$. The quotient ${E_8^p}'/ E_8^p$ is ${\mathbb Z}_p^m$ and ${E_8^p}'(p)=E_8^p$. We give more details on these lattices in the last section.

Let $\Lambda$ be the Leech lattice and $\Lambda^p$ the sublattice fixed by an automorphism of cycle shape $1^mp^m$ where $p$ is a prime such that $m=24/(p+1)$ is an integer. For $p=2$ resp. $p=3$ the lattice is the Barnes-Wall resp. Coxeter-Todd lattice. $\Lambda^p$ has similar properties as the Leech lattice. It has no roots and ${\Lambda^p}'(p)=\Lambda^p$. Furthermore $\Lambda^p$ has level $p$ and ${\Lambda^p}'/\Lambda^p={\mathbb Z}_p^m$. 

Finally we remark that in the following we will adopt to Borcherds' notation in \cite{B3}. For example if $L$ is a lattice that is not positive definite then we define a root in $L$ as a primitive vector of negative rather than positive norm such that the corresponding reflection is an automorphism of $L$. 

\subsection{The singular theta correspondence}

Let $M$ be an even lattice of signature $(b^+,b^-)$ and 
\[ F(\tau)=\sum_{\gamma\in M'/M} f_{\gamma}(\tau)e^{\gamma} \]
be a function on the upper halfplane ${\mathbb H}$ with values in the group ring ${\mathbb C}[M'/M]$. We say that $F$ is holomorphic on ${\mathbb H}$ and meromorphic at the cusps if the components can be written  
\[ f_{\gamma}(\tau) = \sum_{n\in \mathbb Q}c_{\gamma}(n)q^n \]
with $c_{\gamma}(n)=0$ for $n<\!\!<0$. $F$ is a modular form of type $\rho_M$ and weight $m$ if the components satisfy 
\beann 
f_{\gamma}(T\tau) &=&  e(\gamma^2/2)\, f_{\gamma}(\tau) \\
f_{\gamma}(S\tau) &=& \frac{\sqrt{i}^{\,b^--b^+}}{\sqrt{|M'/M|}}\, {\tau^m}
                  \sum_{\delta\in M'/M} e(-(\gamma,\delta))\, f_{\delta}(\tau)
\eeann
under the standard generators 
$S=\left( {0 \atop 1} {-1 \atop 0} \right)$ and 
$T=\left( {1 \atop 0} {1 \atop 1} \right) $ 
of $SL_2({\mathbb Z})$. 

In the next section we will construct vector valued modular forms using products of Dedekind's eta function $\eta(\tau)=q^{1/24}\prod_{n>0}(1-q^n)$ as components. The following lemma will be used to calculate their $S$-transformations. 

\begin{l1} \label{eta}
Let $f(\tau)=\eta\big((k\tau+j)/m\big)$ where $j,k$ and $m$ are integers and let $j',k'$ and $m'$ be integers such that the matrix 
\[ A = \left( \begin{array}{ll} j/k' & -(jj'+kk')/km \\
                                m/k' & -j'/k         \end{array}\right) \]
is in $SL_2({\mathbb Z})$ and $m/k'>0$. Then 

$f(S\tau)=\varepsilon(A) \sqrt{m\tau/m'i} \, \eta\big((k'\tau+j')/m'\big)$ with $\varepsilon(A)$ as given in 
eq. (74.93) of \cite{R}. 
\end{l1}
{\em Proof:} Let 
$F_m=\left( {0 \atop m} {-1 \atop 0} \right)$ be the Fricke involution. Then 
$\eta\big((k\tau+j)/m\big)=\eta(F_mST^jSF_k\tau)$. The lemma now follows from 
$F_mST^jSF_kS=AF_{m'}ST^{j'}SF_{k'}$ and the transformation formula of the eta function.

Borcherds' singular theta correspondence constructs an automorphic form $\Phi_M$ from the vector valued modular form $F$. The singularities of $\Phi_M$ lie on points of the form $\lambda^{\bot}$ where $\lambda$ is a negative norm vector in $M'$ with nonzero coefficient $c_{\lambda}(\lambda^2/2)$. 

If $M$ has signature $(2,b^-)$ then there is an explicit product expansion for $exp(\Phi_M)$ (Theorem 13.3 in \cite{B3}).
\begin{thl1} \label{st1}
Let M be an even lattice of signature $(2,b^-)$ and $F$ a modular form of weight $1-b^-/2$ and representation $\rho_M$ which is holomorphic on ${\mathbb H}$ and meromorphic at the cusps and whose coefficients $c_\lambda(m)$ are integers for $m\leq 0$. Then there is a meromorphic function $\Psi_M(Z_M,F)$ for $Z\in P$ with the following properties.
\begin{enumerate}
\item $\Psi_M(Z_M,F)$ is an automorphic form of weight $c_0(0)/2$ for the group \linebreak 
$Aut(M,F)$ with respect to some unitary character character.
\item The only zeros or poles of $\Psi_M$ lie on the rational quadratic divisors $\lambda^{\bot}$ for $\lambda\in M$ with $\lambda^2<0$ and are zeros of 
order 
\[ \sum_{0<x  \atop x\lambda \in M'}c_{x\lambda}(x^2\lambda^2) \]
or poles if this number is negative.
\item $\Psi_M$ is a holomorphic function if the orders of all zeros are nonnegative. If in addition M has dimension at least 5, or if M has dimension 4 and contains no 2 dimensional isotropic sublattice, then $\Psi_M$ then is a holomorphic automorphic form. If in addition $c_0(0)=b^--2$ then $\Psi_M$ has singular weight and the only nonzero Fourier coefficients of $\Psi_M$ correspond to norm $0$ vectors in $L=K/{\mathbb Z}z$ with $K=M \cap z^{\bot}$.
\item For each primitive norm $0$ vector z in M and for each Weyl chamber W of $L$ the restriction $\Psi_z(Z,F)$ has an infinite product expansion converging when Z is in the neighborhood of the cusp of z and $Y\in W$ which is up to a constant
\[ e((Z,\rho(L,W,F_L))) \prod_{\lambda \in L' \atop (\lambda,W) >0}
               \prod_{\delta \in M'/M \atop \delta|K=\lambda}
      \big( 1-e((\lambda,Z)+(\delta,z')) \big)^{c_{\delta}(\lambda^2/2)}\, . \]
\end{enumerate}
\end{thl1}
If $M$ is Lorentzian the singular theta correspondence gives some information about the automorphism group of $M$ (Theorem 12.1 in \cite{B3}).
\begin{thl1} \label{st2}
Let $M$ be a Lorentzian lattice of dimension $1+b^-$ and $F$ a modular form of weight $(1-b^-)/2$ and representation $\rho_M$ which is holomorphic on ${\mathbb H}$ and meromorphic at the cusps and whose coefficients $c_\lambda(m)$ are real for $m<0$. Suppose that if $\lambda$ is a negative norm vector in $M'$ and $c_\lambda(\lambda^2/2)\neq 0$ then reflection in $\lambda^\bot$ is in $Aut(M,F,C)$. Then $Aut(M,F,C)$ is the semidirect product of a reflection subgroup and a subgroup fixing the Weyl vector $\rho(M,W,F)$ of a Weyl chamber $W$. In particular if the Weyl vector has positive norm then the reflection group of $M$ has finite index in the automorphism group and has only finitely many simple roots. If the Weyl vector has 0 norm but is nonzero then the quotient of the automorphism group of $M$ by the reflection group has a free abelian subgroup of finite index. 
\end{thl1}

\section{\sloppy{Automorphic forms and the fake monster \mbox{algebra}}}

In this section we recall some results about the fake monster algebra and the twisted denominator identities. Then we give a new proof of the twisted denominator identities corresponding to certain automorphisms of prime order and show that they define automorphic forms of singular weight. The idea of the proof is to find appropriate lattices and modular forms and then apply the singular theta correspondence. 
 
\subsection{The fake monster algebra}

The fake monster algebra is a Lie algebra constructed by Borcherds describing the physical states of a bosonic string moving on a torus. It was the first explicit example of a generalized Kac-Moody algebra. We sketch two different constructions. The vertex algebra of the even unimodular Lorentzian lattice $II_{1,25}$ carries an action of the Virasoro algebra. Let $G$ be the space of vectors which are annihilated by the $L_n$ with $n>0$ and have $L_0$-eigenvalue $1$ divided by the kernel of a natural bilinear form. Then the vertex algebra induces a product on $G$ turning it into a Lie algebra called the fake monster algebra. The other construction of this algebra uses the BRST operator. The vertex algebra of the integral lattice $II_{1,25}\oplus {\mathbb Z}^2$ is acted on by the BRST operator $Q$ satisfying $Q^2=0$. Here the space of physical states is given by the cohomology group $Ker(Q)/Im(Q)$. Again the vertex algebra induces the Lie bracket on this space. This Lie algebra is isomorphic to $G$. 

The fake monster algebra has the following properties. The root lattice is the Lorentzian lattice $II_{1,25}=\Lambda(-1)\oplus II_{1,1}$ with elements $\al=(r,m,n)$ and norm $\al^2=r^2-2mn$. A nonzero vector $\al \in II_{1,25}$ is a root if and only if $\al^2\geq -2$. The multiplicity of a root $\al$ is given by $c(\al^2/2)$ where $c(n)$ is the coefficient of $q^n$ in the $1/\Delta(\tau)=1/\eta(\tau)^{24}=q^{-1}+24+324q+3200q^2+\ldots\, $. $1/\Delta(\tau)$ is a modular form for $SL_2({\mathbb Z})$ of weight $-12$ with singularities at the cusps $0$ and $i\infty$. The real simple roots of the fake monster algebra are the norm $-2$ vectors $\al\in II_{1,25}$ with $(\rho,\al)=-1$ where $\rho=(0,0,1)$ is the Weyl vector and the imaginary simple roots are the positive multiples $n\rho$ of the Weyl vector. The Weyl group $G$ is the reflection group of $II_{1,25}$. The denominator identity now reads 
\[ e^{\rho}\prod_{\al \in II_{1,25}^+} \big(1-e^{\al })^{c(\al^2/2)} 
     = \sum_{w\in G}det(w) w\Big(  
             e^{\rho}\prod_{n>0}\big(1-e^{n\rho}\big)^{24}\Big)\, .   \]
The sum in this identity defines the denominator function of the fake monster algebra. It is an automorphic form for $Aut(II_{2,26})$. 

The no-ghost theorem gives an isomorphism from the above descriptions of the physical states to the light-cone states. The automorphism group of the Leech lattice has a natural action on the light-cone states and therefore also on the fake monster algebra. By applying elements of this group to the denominator identity of the fake monster algebra we can calculate twisted denominator identities \cite{B2}.

We consider now the case that the automorphism has cycle shape $1^mp^m$ where $p$ is a prime and $m=24/(p+1)$ an integer (cf. \cite{B2} and \cite{N}). Then the corresponding twisted denominator identity is
\beann  
\lefteqn{
      e^{\rho}
      \prod_{\al\in L^+} \big(1-e^{\al})^{c(\al^2/2)}
      \prod_{\al\in pL'^+} \big(1-e^{\al})^{c(\al^2/2p)}  } \\
 & \qquad = & \sum_{w\in G} det(w) w\Big( e^{\rho} 
     \prod_{n>0}\big(1-e^{n\rho}\big)^m\big(1-e^{pn\rho}\big)^m \Big)\, . 
\eeann
with $L=\Lambda^p(-1)\oplus II_{1,1}$ and $\sum c(n)q^n =(\eta(\tau)\eta(p\tau))^{-m}$. The Weyl vector is $\rho=(0,0,1)$ and the Weyl group $G$ is the reflection group of $L$. It is generated by norm $-2$ vectors in $L$ and the norm $-2p$ vectors in $pL'\subset L$. This identity is also the untwisted denominator identity of a generalized Kac-Moody algebra. The real simple roots of this algebra are the roots $\al$ satisfying $(\rho,\al)=\al^2/2$ and the imaginary simple roots are the positive multiples $n\rho$ of the the Weyl vector with multiplicity $2m$ if $p$ divides $n$ and $m$ else. The root lattice of the algebra is $L$ and the multiplicity of a root $\al$ is given by $c(\al^2/2)$ if $\al$ is in $L$ but not in $pL'$ and by $c(\al^2/2)+c(\al^2/2p)$ if $\al$ is in $pL'$.

\subsection{Automorphic forms}

Now we give a proof of the twisted denominator identities of the fake monster algebra corresponding to automorphisms of cycle shape $1^mp^m$ using the singular theta correspondence. First we work out the case $p=2$ explicitly. The general case will be a simple generalization of this example.

Let
\[ f(\tau)=(\eta(\tau)\eta(2\tau))^{-8}
          =q^{-1}+8+52q+256q^2+1122q^3+4352q^4+\ldots \, . \]
Then $f$ is a modular form for $\Gamma_0(2)$ of weight $-8$ with singularities at the cusps $0$ and $i\infty$. 
The Fourier expansion of 
\[ f(\tau/2) = (\eta(\tau/2)\eta(\tau))^{-8} 
             =  q^{-1/2} + 8 + 52q^{1/2} + 256q + 1122q^{3/2} +\ldots      \] 
can be decomposed into two series with integral and half-integral exponents in $q$. Define 
\[ g_0(\tau) =\big( f(\tau/2)+f((\tau+1)/2) \big)/2 
             = 8 + 256q + 4352q^2 + \ldots \]
and 
\[ g_1(\tau) =\big( f(\tau/2)-f((\tau+1)/2)\big)/2 
             =q^{-1/2} + 52q^{1/2} + 1122 q^{3/2} + \ldots \, . \] 
We will use the functions $f,g_0$ and $g_1$ to construct the vector valued modular function $F$. Therefore we need their transformation properties under the generators of $SL_2({\mathbb Z})$. $f$ and $g_0$ are invariant under $T$ and $g_1(T\tau)=-g_1(\tau)$.   
Lemma \ref{eta} implies the following $S$-transformations 
\beann
f(S\tau)      &=& 2^4 f(\tau/2)/\tau^8 \\
f((S\tau)/2)  &=& f(\tau)/2^4\tau^8 \\
f((S\tau+1)/2)  &=& f((\tau+1)/2)/\tau^8
\eeann
so that
\beann
  f(S\tau) &=& 2^4 \big(g_0(\tau) + g_1(\tau)\big)/ \tau^8 \\
g_0(S\tau) &=& \big( f(\tau)/2^4 + g_0(\tau) - g_1(\tau) \big)/ 2\tau^8 \\
g_1(S\tau) &=& \big( f(\tau)/2^4 - g_0(\tau) + g_1(\tau) \big)/ 2\tau^8 \, .
\eeann

$\Lambda^2(-1)$ is the Barnes-Wall lattice with all norms multiplied by $-1$. We describe the discriminant form of $\Lambda^2(-1)$ in more detail. (Note that the description of $\Lambda^2(-1)'/\Lambda^2(-1)$ in \cite{N} is false.) $\Lambda^2(-1)'/\Lambda^2(-1)={\mathbb Z}_2^8$ has $135$ nonzero elements of norm $0$ mod $2$ and $120$ elements of norm $1$ mod $2$. If $\gamma\in \Lambda^2(-1)'/\Lambda^2(-1)$ is a nonzero element of norm $0$ mod $2$ then there are $71$ nonzero elements $\delta$ of even norm such that $(\gamma,\delta)=0$ mod $1$ and $56$ elements $\delta$ of odd norm such that $(\gamma,\delta)=0$ mod $1$. If $\gamma$ has norm $1$ then $\gamma$ has scalar product $0$ with $63$ of the $135$ nonzero elements of even norm and with $64$ of the $120$ elements of odd norm. The lattice $M=\Lambda^2(-1) \oplus II_{1,1}(2) \oplus II_{1,1}$ is an even lattice of level 2, determinant $2^{10}$ and signature $(2,18)$. The scalar products in $M'/M={\mathbb Z}_2^{10}$ can be derived easily. For example there are $527$ nonzero elements of even norm and $496$ elements of odd norm. 

Now we can define the vector valued modular form $F$. Let 
\[ F(\tau)=\sum_{\gamma\in M'/M} f_{\gamma}(\tau)e^{\gamma} \]  with
\renewcommand{\arraystretch}{1.3}
\[ \begin{array}{lcll} 
f_{\gamma}(\tau) &=& f(\tau) + g_0(\tau)  & \mbox{if $\gamma = 0$}     \\
                 &=& g_0(\tau) & \mbox{if $\gamma^2/2 =0$ mod $1$} \\
                 &=& g_1(\tau) & \mbox{if $\gamma^2/2 =1/2$ mod $1$.}
\end{array} \]
Then $F$ is a modular form of weight $-8$ and representation $\rho_M$. The $T$-invariance is clear. $F$ also transforms correctly under $S$. We show this for $f_{0}(\tau)$.
\beann 
2^5\tau^8 f_{0}(S\tau) &=& 2^5\tau^8 f(S\tau) + 2^5\tau^8 g_0(S\tau) \\
                &=& f(\tau) + (2^9+2^4)g_0(\tau) + (2^9-2^4)g_0(\tau)\\
                &=& \sum f_{\gamma}(\tau)e^{\gamma}
\eeann
by the above formulas. The proof for the other components is similar. Note that there are only $2$ other cases which must be considered. 

The singular theta correspondence now implies that there is a holomorphic automorphic form $\Psi_M$ for $Aut(M)$ of weight $16/2=8$. The zeros of $\Psi_M$ are zeros of order 1 coming from divisors of norm $-2$ vectors of $M$ and norm $-1$ vectors of $M'$.

The lattice $M$ has $2$ orbits of primitive norm $0$ vectors under $Aut(M)$, 
one of level $1$ and one of level $2$. We explain the terminology below. Near these cusps we can expand $\Psi_M$ in infinite products. Since $\Psi_M$ has singular weight the nonzero Fourier coefficients of $\Psi_M$ correspond to norm $0$ vectors. This allows us to work out the Fourier expansions at the cusps explicitly.

Define numbers $c(n)$ by
\[ \sum c(n)q^n = f(\tau) = q^{-1}+8+52q+256q^2+1122q^3+4352q^4+\ldots \, . \]  
Level $1$ cusp: We decompose $M=L\oplus II_{1,1}$ where $L=\Lambda^2(-1) \oplus II_{1,1}(2)$ and take $z$ as primitive norm $0$ vector in $II_{1,1}$. Then the product expansion of $\Psi_z(Z,F)$ is 
\[ e((\rho,Z))
   \prod_{\lambda\in L^+}  \big(1-e((\lambda,Z))\big)^{c(\lambda^2/2)}
   \prod_{\lambda\in L'^+} \big(1-e((\lambda,Z))\big)^{c(\lambda^2)} \]
\[ =\sum_{w\in G}det(w)e((w\rho,Z)) 
          \prod_{n>0}\big(1-e((n w\rho,Z))\big)^8 
                     \big(1-e((2nw\rho,Z))\big)^8 \]
where $\rho=(0,0,1/2)$ and $G$ is the reflection group generated by norm $-1$ vectors of $L'$ and the norm $-2$ vectors of $L\subset L'$. Note that the vectors of $L'$ have integral norms because $L$ has level $2$. 

Level $2$ cusp: Here we write $M=L \oplus II_{1,1}(2)$ with $L=\Lambda^2(-1)\oplus II_{1,1}$ and take $z$ as primitive norm $0$ vector in $II_{1,1}(2)$. We say that $z$ has level $2$ because $|M'/M|=2^2 |L'/L|$. At this cusp the product expansion of $\Psi_z(Z,F)$ is 
\[ e((\rho,Z)) 
   \prod_{\lambda\in L^+} \big(1-e((\lambda,Z))\big)^{c(\lambda^2/2)}
   \prod_{\lambda\in L'^+} \big(1-e((\lambda,Z))\big)^{c(\lambda^2/4)} \]
\[   = \sum_{w\in G}det(w)
           e(( w \rho,Z)) 
          \prod_{n>0}\big(1-e((n w\rho,Z))\big)^8 
                     \big(1-e((2nw\rho,Z))\big)^8     \]
where $\rho=(0,0,1)$ and $G$ is the reflection group generated by the norm $-2$ vectors of $L$ and the norm $-4$ vectors of $2L'\subset L$. 

The expansion of $\Psi_M$ at the level $2$ cusp is the twisted denominator identity of the fake monster algebra corresponding to the prime $p=2$. This gives a new proof of this identity and shows that the denominator function of the corresponding generalized Kac-Moody algebra is an automorphic form of singular weight. 

The two expansions of $\Psi_M$ look rather similar and are really the same. If we rescale the dual of $\Lambda^2(-1) \oplus II_{1,1}(2)$ with a factor $2$ then the expansion of $\Psi_M$ at the level $1$ cusp goes over into the expansion at the other cusp.

No we turn to the general case. Let $p$ be a prime such that $m=24/(p+1)$ is an integer. The eta product 
\[ f(\tau) = (\eta(\tau)\eta(p\tau))^{-m} = q^{-1} + m + \ldots \]
is a modular form of weight $-m$ for 
$\Gamma_0(p)=\{ \big( {a \atop c} {b \atop d} \big) 
      \in SL_2({\mathbb Z})\, |\, c=0 \mbox{ mod } p \}$ 
if $m$ is even and for 
$\Gamma_1(p)=\{ \big( {a \atop c} {b \atop d} \big) \in SL_2({\mathbb Z})
 \, |\, a=d=1 \mbox{ mod } p, \, c=0 \mbox{ mod } p \}$ 
if $m$ is odd. $f$ has singularities at the cusps $0$ and $i \infty$. 
$f(\tau/p)$ can be written 
\[ f(\tau/p) = g_0(\tau)+g_1(\tau)+\ldots+g_{p-1}(\tau)  \]
where the functions $g_j$ have Fourier expansions of the form $\sum a(n) q^{n+j/p}$ with $n\in {\mathbb Z}$. We will use the functions $f, g_0, \ldots, g_{p-1}$ to construct a vector valued modular function. The $T$-transformations of these functions are clear. As above their $S$-transformations can be calculated with Lemma \ref{eta}. 

Let $M=\Lambda^p(-1) \oplus II_{1,1}(p) \oplus II_{1,1}$. Then $M$ is an even lattice of level $p$, determinant $p^{m+2}$ and signature $(2,2m+2)$. We define 
\[ F(\tau)=\sum_{\gamma\in M'/M} f_{\gamma}(\tau)e^{\gamma} \]  with
\renewcommand{\arraystretch}{1.3}
\[ \begin{array}{lcll} 
f_{\gamma}(\tau) &=& f(\tau) + g_0(\tau)  & \mbox{if $\gamma = 0$}     \\
                 &=& g_j(\tau) & \mbox{if $\gamma^2/2 =j/p$ mod $1$}    
\end{array} \]
Then we have
\begin{p1}
$F$ is a modular form of weight $-m$ and representation $\rho_M$ which is holomorphic on ${\mathbb H}$ and meromorphic at the cusps. 
\end{p1}
{\em Proof:}
$F$ clearly transforms correctly under $T$. We can work out the quotient $M'/M$ by using an explicit description of the Leech lattice. A case by case analysis then shows that $F$ also transforms correctly under $S$. This proves the proposition.

Now the singular theta correspondence implies 
\begin{thp1}
There is a holomorphic automorphic form $\Psi_M$ for $Aut(M)$ of weight $m$. The zeros of $\Psi_M$ are zeros of order $1$ coming from divisors of norm $-2$ vectors of $M$ and norm $-2/p$ vectors of $M'$. $\Psi_M$ has singular weight so that the only nonzero Fourier coefficients of $\Psi_M$ correspond to norm $0$ vectors. 

Let $c(n)$ be given by $f(\tau)=(\eta(\tau)\eta(p\tau))^{-m}=\sum c(n) q^n$. Then $\Psi_M$ has the following expansions at the cusps. 

At the level $1$ cusp we decompose $M=L \oplus II_{1,1}$ with $L=\Lambda^p(-1)+II_{1,1}(p)$ and take $z$ as primitive norm $0$ vector in $II_{1,1}$. Then the product expansion of $\Psi_z(Z,F)$ is 
\[ e((\rho,Z))
   \prod_{\lambda\in L} \big(1-e((\lambda,Z))\big)^{c(\lambda^2/2)}
   \prod_{\lambda\in L'} \big(1-e((\lambda,Z))\big)^{c(p\lambda^2/2)} \]
\[ =\sum_{w\in G}det(w)e((w\rho,Z)) 
          \prod_{n>0}\big(1-e((n w\rho,Z))\big)^m 
                     \big(1-e((pnw\rho,Z))\big)^m  \] 
where $\rho=(0,0,1/p)$ and $G$ is the reflection group generated by norm $-2/p$ vectors of $L'$ and the norm $-2$ vectors of $L$. 

At the level $p$ cusp we write $M=L\oplus II_{1,1}(p)$ with $L=\Lambda^p(-1)\oplus II_{1,1}$ so that $|M'/M|=p^2 |L'/L|$ and take $z$ as primitive norm $0$ vector in $II_{1,1}(p)$. Here the product expansion of $\Psi_z(Z,F)$ is 
\[ e((\rho,Z)) 
   \prod_{\lambda\in L} \big(1-e((\lambda,Z))\big)^{c(\lambda^2/2)}
   \prod_{\lambda\in pL'} \big(1-e((\lambda,Z))\big)^{c(\lambda^2/2p)} \]
\[   = \sum_{w\in G}det(w)e((w\rho,Z)) 
          \prod_{n>0}\big(1-e((n w\rho,Z))\big)^m 
                     \big(1-e((pnw\rho,Z))\big)^m   \] 
where $\rho=(0,0,1)$ and $G$ is the reflection group generated by the norm $-2$ vectors of $L$ and the norm $-2p$ vectors of $pL'\subset L$.

The two expansions are identical upon rescaling the lattice $(\Lambda^p(-1)+II_{1,1}(p))'$ by $p$. 
\end{thp1}
{\em Proof:} The only thing left is to calculate the Fourier expansions of $\Psi_M$. We can do this using the fact that the nonzero Fourier coefficients correspond to norm $0$ vectors of $L$. This proves the theorem.
\begin{cp1}
The denominator function of the generalized Kac-Moody algebra obtained by twisting the fake monster algebra with an automorphism of cycle shape $1^mp^m$ defines a holomorphic automorphic form of singular weight.
\end{cp1}
We remark that the results in this section can be easily extended to the case $p=1$.

\section{\sloppy{Automorphic forms and the fake monster \mbox{superalgebra}}}

In this section we prove analogous results as in section $3$ for the fake monster superalgebra. The main difference is that the expansions of the automorphic forms at the $2$ different cusps do not coincide so that we get new infinite product identities which are denominator identities of generalized Kac-Moody superalgebras.

\subsection{The fake monster superalgebra}

The fake monster superalgebra \cite{S1} is a supersymmetric generalized Kac-Moody superalgebra describing the physical states of a superstring moving on a torus. It can be constructed similar to the fake monster algebra as the cohomology group of a BRST operator acting on the vertex algebra of a rational 18 dimensional lattice. The fake monster superalgebra has root lattice $II_{1,9}=E_8(-1)\oplus II_{1,1}$ with elements $\al=(r,m,n)$ and norm $\al^2=r^2-2mn$. The roots are the nonzero vectors $\al$ with $\al^2\geq 0$. The multiplicity of a root $\al$ is given by $mult_0(\al)=mult_1(\al)=c(-\al^2/2)$ where $c(n)$ is the coefficient of $q^n$ in 
$ 8 \eta(2\tau)^8/\eta(\tau)^{16} = 8+128q+1152q^2+\ldots$ which is a modular form for $\Gamma_{0}(2)$ of weight $-8$. The simple roots of the fake monster superalgebra are the norm $0$ vectors in the closure of the positive cone of $II_{1,9}$. This implies that the Weyl group is trivial and the Weyl vector is $0$. The denominator identity is given by 
\[ \prod_{\al\in II_{1,9}^{+}}
     \frac{ (1-e^{\al})^{c(\al^2/2)} }
          { (1+e^{\al})^{c(\al^2/2)} }
   =  1 + \sum a(\lambda)e^{\lambda}  \]
where $a(\lambda)$ is the coefficient of $q^n$ in 
\[  \frac{\eta(\tau)^{16}}{\eta(2\tau)^8}
 =  1 - 16q + 112q^2 - 448q^3 + 1136q^4  - \ldots  \]
if $\lambda$ is $n$ times a primitive norm $0$ vector in $II_{1,9}^{+}$ and $0$ else. 

Using the no-ghost theorem we can construct an action of $2.Aut(E_8)$ on the fake monster superalgebra and calculate twisted denominator identities \cite{S2}. The identity corresponding to an automorphism of cycle shape $1^mp^m$ with  $p$ prime and $m=8/(p+1)$ integral is 
\[ \prod_{\al\in L^{+}}
      \frac{ (1-e^{\al})^{c(\al^2/2)} }
           { (1+e^{\al})^{c(\al^2/2)} }
   \prod_{\al \in pL'^{+}}
      \frac{ (1-e^{\al})^{c(\al^2/2p)} }
           { (1+e^{\al})^{c(\al^2/2p)} }
=  1 + \sum a(\lambda)e^{\lambda}  \]
where $L=E_8^p(-1)\oplus II_{1,1}$ and 
$\sum c(n)q^n = m(\eta(2p\tau)\eta(2\tau))^m/ (\eta(p\tau)  \eta(\tau))^{2m}$.
$a(\lambda)$ is the coefficient of $q^n$ in 
$ ( \eta(p\tau) \eta(\tau)  )^{2m} /( \eta(2p\tau)\eta(2\tau) )^m$
if $\lambda$ is $n$ times a primitive norm $0$ vector in $L^+$ and $0$ else. This is the untwisted denominator identity of a supersymmetric generalized Kac-Moody superalgebra whose simple roots are the norm $0$ vectors in the closure of the positive cone of $L$. The even and the odd multiplicity of a simple root $\lambda$ is $2m$ if $\lambda$ is $n$ times a primitive vector with $p$ dividing $n$ and $m$ else.

\subsection{Automorphic forms}

Now we construct the automorphic forms whose expansions give the denominator identities of the twisted fake monster superalgebras and some new identities. In contrast to the bosonic case the components of the vector valued modular form do not only depend on the norm of the elements in $M'/M$ but also on their order. 

We will start with the case corresponding to $p=3$ as an example. Let 
\[  f(\tau) = 
    2\, \frac{\eta(6\tau)^2\eta(2\tau)^2}{\eta(3\tau)^4\eta(\tau)^4} =
    2 + 8q + 24q^2+ 72q^3 + 184q^4 + \ldots \]
and 
\[ \gamma(\tau) = 
        \frac{\eta(3\tau)^2\eta(\tau/2)^2}{\eta(3\tau)^4\eta(\tau)^4} =
        q^{-1/2} - 2 + 3q^{1/2} - 8q + 15q^{3/2} - 24q^2 + \ldots  \, .\]
$f$ is a modular form for $\Gamma_0(6)$ of weight $-2$. 

$f$ and $\gamma$ are related by supersymmetry \cite{S2} which means that the Fourier expansion of 
\[ \delta(\tau) = f(\tau) + \gamma(\tau) = q^{-1/2} + 3q^{1/2} + 15q^{3/2} 
                                           + 43q^{5/2} + \ldots  \] 
only contains half-integral powers of $q$.
We write 
\[ f(\tau/3) =  g_0(\tau) + g_2(\tau) + g_4(\tau) \]
with 
\beann   
   g_0(\tau) & = & \big( f(\tau/3)+f((\tau+1)/3)+f((\tau+2)/3)\big)/3 \\
             & = & 2 + 72q + 984 q^2 + \ldots \\
   g_2(\tau) & = & \big( f(\tau/3)+\varepsilon^2 f((\tau+1)/3)
                                +\varepsilon f((\tau+2)/3)\big)/3 \\
             & = & 8q^{1/3} + 184 q^{4/3} + 2112 q^{7/3} + \ldots \\
   g_4(\tau) & = & \big( f(\tau/3)+\varepsilon f((\tau+1)/3)
                                +\varepsilon^2 f((\tau+2)/3)\big)/3 \\
             & = & 24q^{2/3} + 432 q^{5/3} + 4344 q^{8/3} + \ldots
\eeann
and 
\[ \delta(\tau/3) =  \lambda_1(\tau) +  \lambda_3(\tau) +  \lambda_5(\tau) \]
with 
\beann
   \lambda_1(\tau) &=& 
      \big( \delta(\tau/3)+\varepsilon^2\delta((\tau+2)/3)
                             +\varepsilon\delta((\tau+4)/3)\big)/3 \\
                   &=& 3 q^{1/6} + 114 q^{7/6} + 1437 q^{13/6} + \ldots \\
   \lambda_3(\tau) &=& 
      \big( \delta(\tau/3)+\delta((\tau+2)/3)+\delta((\tau+4)/3)\big)/3 \\
                   &=& 15 q^{3/6} + 285q^{9/6} + 3051 q^{15/6}  + \ldots \\ 
   \lambda_5(\tau) &=& 
       \big( \delta(\tau/3)+\varepsilon\delta((\tau+2)/3)
                             +\varepsilon^2\delta((\tau+4)/3)\big)/3 \\
                   &=& q^{-1/6} + 43 q^{5/6} + 662 q^{11/6} + \ldots \, .
\eeann 
Note that $f((\tau+4)/3)=f((\tau+1)/3)$. We will use these functions to construct the vector valued modular function $F$. We determine the transformation properties under the generators of $SL_2({\mathbb Z})$. The transformations under $T$ are clear. Using Lemma \ref{eta} we find for the $S$-transformations 
\beann 
f(S\tau)           &=& -3\gamma(\tau/3)/2\tau^2       \\
f((S\tau)/3)       &=& -\gamma(\tau)/6\tau^2          \\
f((S\tau+2)/3)     &=& \gamma((\tau+4)/3)/2\tau^2    \\
f((S\tau+4)/3)     &=& \gamma((\tau+2)/3)/2\tau^2    
\eeann
and 
\beann
\gamma(S\tau)        &=& -6f(\tau/3)/\tau^2     \\
\gamma((S\tau)/3)    &=& -2f(\tau)/3\tau^2      \\
\gamma((S\tau+2)/3)  &=& 2f((\tau+4)/3)/\tau^2   \\
\gamma((S\tau+4)/3)  &=& 2f((\tau+2)/3)/\tau^2   
\eeann
so that 
\beann
f(S\tau)   &=&  3(g_0(\tau)+g_2(\tau)+g_4(\tau))/2\tau^2 
               -3(\lambda_1(\tau)+\lambda_3(\tau)+\lambda_5(\tau))/2\tau^2  \\
g_0(S\tau) &=& (f(\tau)-\delta(\tau))/18\tau^2 
               +(\lambda_3(\tau)-g_0(\tau))/2\tau^2 \\
           & & +((g_0(\tau)+g_2(\tau)+g_4(\tau))/6\tau^2
               -(\lambda_1(\tau)+\lambda_3(\tau)+\lambda_5(\tau))/6\tau^2\\
g_2(S\tau) &=& (f(\tau)-\delta(\tau))/18\tau^2 
               +(\lambda_1(\tau)-g_4(\tau))/2\tau^2 \\
           & & +((g_0(\tau)+g_2(\tau)+g_4(\tau))/6\tau^2
               -(\lambda_1(\tau)+\lambda_3(\tau)+\lambda_5(\tau))/6\tau^2\\
g_4(S\tau) &=& (f(\tau)-\delta(\tau))/18\tau^2 
               +(\lambda_5(\tau)-g_2(\tau))/2\tau^2 \\
           & & +((g_0(\tau)+g_2(\tau)+g_4(\tau))/6\tau^2
               -(\lambda_1(\tau)+\lambda_3(\tau)+\lambda_5(\tau))/6\tau^2
\eeann
and
\beann
\delta(S\tau) &=& -9(g_0(\tau)+g_2(\tau)+g_4(\tau))/2\tau^2
                  -3(\lambda_1(\tau)+\lambda_3(\tau)+\lambda_5(\tau))/2\tau^2\\
\lambda_1(S\tau) &=& -(3f(\tau)+\delta(\tau))/18\tau^2
                     +(3g_2(\tau)+\lambda_5(\tau))/2\tau^2 \\
            & & -(g_0(\tau)+g_2(\tau)+g_4(\tau))/2\tau^2
                -(\lambda_1(\tau)+\lambda_3(\tau)+\lambda_5(\tau))/6\tau^2 \\
\lambda_3(S\tau) &=& -(3f(\tau)+\delta(\tau))/18\tau^2
                     +(3g_0(\tau)+\lambda_3(\tau))/2\tau^2 \\
            & & -(g_0(\tau)+g_2(\tau)+g_4(\tau))/2\tau^2
                -(\lambda_1(\tau)+\lambda_3(\tau)+\lambda_5(\tau))/6\tau^2 \\
\lambda_5(S\tau) &=& -(3f(\tau)+\delta(\tau))/18\tau^2
                     +(3g_4(\tau)+\lambda_1(\tau))/2\tau^2 \\
            & & -(g_0(\tau)+g_2(\tau)+g_4(\tau))/2\tau^2
                -(\lambda_1(\tau)+\lambda_3(\tau)+\lambda_5(\tau))/6\tau^2 
\eeann

We describe the discriminant form of $E_8^3(-1)$. The quotient $E_8^3(-1)'/E_8^3(-1)$ is a 2 dimensional vector space over ${\mathbb Z}_3$. We can choose a basis $\{\gamma_1, \gamma_2\}$ with $\gamma_1^2/2=\gamma_2^2/2=-1/3$ mod $1$ and $(\gamma_1,\gamma_2)=0$ mod $1$.

We define the even lattice $M=E_8^3(-1) \oplus II_{1,1}(6) \oplus II_{1,1}$ of level $6$, signature $(2,6)$ and determinant $324$. The scalar product of the elements in $M'/M= {\mathbb Z}_3^2 \times {\mathbb Z}_6^2$ can be worked out using the above results on $E_8^3(-1)$ and $II_{1,1}(6)$. Now let 
\[ F(\tau)=\sum_{\gamma\in M'/M} f_{\gamma}(\tau)e^{\gamma} \]
with
\renewcommand{\arraystretch}{1.2}
\[ \begin{array}{lcll} 
f_{\gamma}(\tau) &=& f(\tau) + g_0(\tau)  & \mbox{if $\gamma = 0$}     \\
                 &=& -f(\tau) - g_0(\tau) & \mbox{if $\gamma^2/2 = 0$ 
                                                and $\gamma$ has order $2$}  \\
                 &=& \delta(\tau) + \lambda_3(\tau)  
                                          & \mbox{if $\gamma^2/2 = 1/2$ 
                                            and $\gamma$ has order $2$}   \\
                 &=& \lambda_j(\tau)      & \mbox{if $\gamma^2/2 = j/6$ 
                                            where $j$ is odd and 
                                            $\gamma$ has order $6$} \\ 
                 &=&  g_j(\tau)           & \mbox{if $\gamma^2/2 = j/6$ 
                                            where $j$ is even and 
                                            $\gamma$ has order $3$}\\
                 &=& -g_j(\tau)           & \mbox{if $\gamma^2/2 = j/6$ 
                                            where $j$ is even and 
                                            $\gamma$ has order $6$}
\end{array} \]
where as usual $\gamma^2/2$ is taken mod $1$. 
Then $F$ is modular form of weight $-2$ and representation $\rho_M$. Again the $T$-transformations are clear. We show that $f_{0}(\tau)$ behaves correctly under $S$. The following table shows how often the corresponding functions appear as component in $F$. 
\[ \begin{array}{lc}
 f(\tau) + g_0(\tau), \, \delta(\tau) + \lambda_3(\tau) \qquad    & 1  \\
-f(\tau) - g_0(\tau)                                              & 2  \\
 g_0(\tau), \, \lambda_3(\tau)                                    & 20 \\
 g_2(\tau), \, g_4(\tau), \, \lambda_1(\tau), \, \lambda_5(\tau)  & 30 \\
-g_0(\tau)                                                        & 40 \\
-g_2(\tau), \, -g_4(\tau)                                         & 60
\end{array} \]
We have 
\beann
-18 \tau^2 f_0(S\tau) &=& -18 \tau^2 f(S\tau) - 18 \tau^2 g_0(S\tau) \\
                      &=&  -f(\tau) 
                          -21g_0(\tau) -30g_2(\tau) -30g_4(\tau) \\
                      & & +\delta(\tau)
                          +30\lambda_1(\tau) 
                          +21\lambda_3(\tau) 
                          +30\lambda_5(\tau) \\
                      &=& \sum f_{\gamma}(\tau)e^{\gamma} \, .
\eeann
The calculations for the other components are similar.

We should remark that the supersymmetry relation between $f(\tau)$ and $\gamma(\tau)$ is essential for $F$ to be a modular form with representation $\rho_M$. 

By the singular theta correspondence there is a holomorphic automorphic form $\Psi_M$ for $Aut(M)$ of weight $4/2=2$. The zeros of $\Psi_M$ are zeros of order 1 coming from divisors of norm $-1/3$ vectors in $M'$ and from divisors $\al$ in $M'$ of norm $-1$ with $2\al \in M$.

The lattice $M$ has primitive norm $0$ vectors of level $1$ and level $6$. We can work out the Fourier expansions at the corresponding cusps using that $\Psi_M$ has singular weight.

Define $c(n)$ by
\[ \sum c(n)q^n = f(\tau)+\delta(\tau) 
   = q^{-1/2} + 2 + 3q^{1/2} + 8q + 15q^{3/2} + 24q^2 
     + \ldots \, . \]

Level $1$ cusp: We decompose $M=L\oplus II_{1,1}$ where $L=E_8^3(-1) \oplus II_{1,1}(6)$ and take $z$ as primitive norm $0$ vector in $II_{1,1}$. Then the product expansion of $\Psi_z(Z,F)$ is 
\[ e((\rho,Z))
    \prod_{\al\in L'^{+} \atop 2\alpha \in L} 
                \big(1-e((\al,Z))\big)^{\pm c(\al^2/2)}
    \prod_{\al\in L'^{+} } 
                \big(1-e((\al,Z))\big)^{\pm c(3 \al^2/2)}  \]
\[ =\sum_{w\in G}det(w)e((w\rho,Z))
          \prod_{n>0}\big(1-e((n w\rho,Z))\big)^{(-1)^n2} 
                     \big(1-e((3nw\rho,Z))\big)^{(-1)^n2} \]
where the sign in the exponent of the first product is $-$ if $\al^2$ is even and the image of $\al$ has even order in $L'/L$, i.e. order $2$, and $+$ else and the sign in the exponent of the second product is $-$ if $3\al^2/2$ is integral and $\al$ has even order and $+$ in the other cases. We remark that $\al\in L'$ with $2\al \in L$ implies $\al^2$ integral because $L$ is even and has level $6$. The Weyl vector is $\rho=(0,0,1/6)$ and the Weyl group $G$ is generated by the $\al$ in $L'$ with $\al^2=-1/3$ and $\al$ in $L'$ of norm $-1$ with $2\al \in L$. We remark that the roots of $L'$ are the vectors $\al$ of norm $-1/3, -2/3, -1$ and $-2$ with resp. $6\al\in L, 3\al\in L, 2\al\in L$ and $\al\in L$. This implies that $G$ has infinite index in the full reflection group of $L'$.

This identity is the denominator identity of a generalized Kac-Moody superalgebra with the following simple roots. The real simple roots are the simple roots of the reflection group $G$, i.e. the roots $\al$ satisfying $(\rho,\al)=\al^2/2$. The imaginary simple roots are the positive multiples $n\rho$ of the Weyl vector with multiplicity $(-1)^n4$ if $3$ divides $n$ and $(-1)^n2$ else. Here we use the convention that odd roots have negative multiplicity. Note that there are no odd real simple roots. The root lattice of this algebra is $L'$. A vector $\al$ in $L'$ is a root if and only if $\al^2\geq -1$. A root $\al$ is odd if and only if $3\al^2/2$ is integral and $\al$ has even order. The multiplicity of an even root $\al\in L'$ is $mult(\al) = c(3\al^2/2)$ if $2\al\notin L$ and $mult(\al) = c(\al^2/2) + c(3\al^2/2)$ if $2\al \in L$. Up to a sign the same formula holds for the odd roots. 

Level $6$ cusp: Here we write $M=L\oplus II_{1,1}(6)$ with $L=E_8^3(-1) \oplus II_{1,1}$ and take $z$ as primitive norm $0$ vector in $II_{1,1}(6)$. Then $|M'/M|=6^2 |L'/L|$ so that $z$ has level $6$. The product expansion of $\Psi_z(Z,F)$ is 
\[ \prod_{\al\in L^{+}}
            \frac{ \big(1-e((\al,Z))\big)^{c(\al^2/2)} }
                 { \big(1+e((\al,Z))\big)^{c(\al^2/2)} } \,
   \prod_{\al \in 3L'^{+}}
            \frac{ \big(1-e((\al,Z))\big)^{c(\al^2/6)} }
                 { \big(1+e((\al,Z))\big)^{c(\al^2/6)} } \]
\[ =  1 + \sum a(\lambda)e((\lambda,Z))  \]
where $a(\lambda)$ is the coefficient of $q^n$ in 
\[  \frac{ \big( \eta(3\tau) \eta(\tau)  \big)^{4} }
         { \big( \eta(6\tau)\eta(2\tau) \big)^2    }
      = 1 - 4q + 4q^2 - 4q^3 + 20q^4 - 24q^5 + 4q^6 - \ldots  \]
if $\lambda$ is $n$ times a primitive norm $0$ vector in $L^+$ and $0$ else.

This is the twisted denominator identity of the fake monster superalgebra corresponding to the prime $p=3$. 

The general case is as follows. Let $p$ be a prime such that $8/(p+1)$ is an integer. Define 
\[  f(\tau) = 
    m\, \frac{\big(\eta(2p\tau)\eta(2\tau)\big)^m}
             {\big(\eta(p\tau)\eta(\tau)\big)^{2m}} =
    m + 2m^2 q + \ldots \]
and 
\[ \gamma(\tau) = 
        \frac{\big(\eta(p\tau/2)\eta(\tau/2)\big)^m}
             {\big(\eta(p\tau)\eta(\tau)\big)^{2m}} =
        q^{-1/2} - m  + \ldots  \, .\]
Then $f$ is a modular form of weight $-m$ for $\Gamma_0(2p)$ if $m$ is even and for $\Gamma_1(2p)$ if $m$ is odd.

A supersymmetry relation implies that $\delta(\tau)=f(\tau)+\gamma(\tau)$ contains only half-integral powers of $q$. 

We decompose $f(\tau/p)$ as 
\[ f(\tau/p) = g_0(\tau)+g_2(\tau)+\ldots+g_{2p-2}(\tau)  \]
such that the functions $g_j$ have Fourier expansions of the form $\sum a(n)q^{n+j/2p}$ with $n\in {\mathbb Z}$. Similarly we write $\delta(\tau/p)$ as
\[ \delta(\tau/p) = \lambda_1(\tau)+\lambda_3(\tau)+\ldots
                    +\lambda_{2p-1}(\tau)  \]
where the $\lambda_j$ have Fourier expansions of the form $\sum a(n)q^{n+j/2p}$ with $n\in {\mathbb Z}$. We will use the functions $f, g_0, \ldots, g_{2p-2}$ and $\delta, \lambda_1, \ldots, \lambda_{2p-1}$ to construct a vector valued modular form. Their transformations under $T$ are clear and their $S$-transformations follow from Lemma \ref{eta}. 

The quotient $E_8^p(-1)'/E_8^p(-1)={\mathbb Z}_p^m$ is a vector space over ${\mathbb Z}_p$ with an orthogonal basis $\{\gamma_1,\ldots,\gamma_m\}$ satisfying $\gamma_j^2/2=-1/p$ mod $1$. Define the even lattice $M=L\oplus II_{1,1}(2p)\oplus II_{1,1}$ of level $2p$, determinant $2^2p^{m+2}$ and signature $(2,2m+2)$. 
Let 
\[ F(\tau)=\sum_{\gamma\in M'/M} f_{\gamma}(\tau)e^{\gamma} \]
with
\renewcommand{\arraystretch}{1.2}
\[ \begin{array}{lcll} 
f_{\gamma}(\tau) &=& f(\tau) + g_0(\tau)  & \mbox{if $\gamma = 0$}     \\
                 &=& -f(\tau) - g_0(\tau) & \mbox{if $\gamma^2/2 = 0$ 
                                                and $\gamma$ has order $2$}  \\
                 &=& \delta(\tau) + \lambda_p(\tau)  
                                          & \mbox{if $\gamma^2/2 = 1/2$ 
                                            and $\gamma$ has order $2$}   \\
                 &=& \lambda_j(\tau)      & \mbox{if $\gamma^2/2 = j/2p$ 
                                            where $j$ is odd and 
                                            $\gamma$ has order $2p$} \\ 
                 &=&  g_j(\tau)           & \mbox{if $\gamma^2/2 = j/2p$ 
                                            where $j$ is even and 
                                            $\gamma$ has order $p$}\\
                 &=& -g_j(\tau)           & \mbox{if $\gamma^2/2 = j/2p$ 
                                            where $j$ is even and 
                                            $\gamma$ has order $2p$}
\end{array} \]
where $\gamma^2/2$ is taken mod $1$. Then we have
\begin{p1}
$F$ is a modular form of weight $-m$ and representation $\rho_M$ which is holomorphic on ${\mathbb H}$ and meromorphic at the cusps. 
\end{p1}
{\em Proof:}
The components $f_\gamma$ are $T$-invariant by construction. The 
$S$-transforma\-tions follow from those of the functions $f,g_j$ and $\delta,\lambda_j$ and the scalar products in $M'/M$ which can be worked out using the descriptions of the discriminant forms of $E_8^p(-1)$ and $II_{1,1}(2p)$. This proves the proposition. 

From the singular theta correspondence we get
\begin{thp1}
There is a holomorphic automorphic form $\Psi_M$ for $Aut(M)$ of weight $m$. The zeros of $\Psi_M$ are zeros of order 1 coming from divisors of norm $-1/p$ vectors in $M'$ and from divisors $\al$ in $M'$ with norm $-1$ and $2\al\in M$. $\Psi_M$ has singular weight so that the only nonzero Fourier coefficients of $\Psi_M$ correspond to norm $0$ vectors. 

Define $c(n)$ by $\sum c(n) q^n=f(\tau)+\delta(\tau)$. Then $\Psi_M$ has the following expansions at the cusps. 

At the level $1$ cusp we decompose $M=L \oplus II_{1,1}$ with $L=E_8^p(-1)\oplus II_{1,1}(2p)$ and take $z$ as primitive norm $0$ vector in $II_{1,1}$. Then the product expansion of $\Psi_z(Z,F)$ is 
\[ e((\rho,Z))
    \prod_{\al\in L'^{+} \atop 2\alpha \in L} 
                \big(1-e((\al,Z))\big)^{\pm c(\al^2/2)}
    \prod_{\al\in L'^{+} } 
                \big(1-e((\al,Z))\big)^{\pm c(p \al^2/2)}  \]
\[ =\sum_{w\in G}det(w)e((w\rho,Z))
          \prod_{n>0}\big(1-e((n w\rho,Z))\big)^{(-1)^nm} 
                     \big(1-e((pnw\rho,Z))\big)^{(-1)^nm} \]
where the sign in the exponent of the first product is $-$ if $\al^2$ is even and the image of $\al$ has even order in $L'/L$, i.e. order $2$, and $+$ else and the sign in the exponent of the second product is $-$ if $p\al^2/2$ is integral and $\al$ has even order and $+$ in the other cases. The Weyl vector is $\rho=(0,0,1/2p)$ and the Weyl group $G$ is generated by the $\al$ in $L'$ with $\al^2=-1/p$ and $\al$ in $L'$ of norm $-1$ with $2\al \in L$.

At the level $2p$ cusp we write $M=L \oplus II_{1,1}(2p)$ with $L = E_8^p(-1) \oplus II_{1,1}$ so that $|M'/M|=2^2p^2 |L'/L|$ and take $z$ as primitive norm $0$ vector in $II_{1,1}(2p)$. Here the product expansion of $\Psi_z(Z,F)$ is
\[ \prod_{\al\in L^{+}}
            \frac{ \big(1-e((\al,Z))\big)^{c(\al^2/2)} }
                 { \big(1+e((\al,Z))\big)^{c(\al^2/2)} } \,
   \prod_{\al \in pL'^{+}}
            \frac{ \big(1-e((\al,Z))\big)^{c(\al^2/2p)} }
                 { \big(1+e((\al,Z))\big)^{c(\al^2/2p)} } \]
\[ =  1 + \sum a(\lambda)e((\lambda,Z))  \]
where $a(\lambda)$ is the coefficient of $q^n$ in 
\[  \frac{ \big( \eta(p\tau) \eta(\tau)  \big)^{2m} }
         { \big( \eta(2p\tau)\eta(2\tau) \big)^m    }
      = 1 - 2mq + \ldots  \]
if $\lambda$ is $n$ times a primitive norm $0$ vector in $L^+$ and $0$ else.
\end{thp1}
{\em Proof:} We only need to calculate the Fourier expansions of the products which can be done using that $\Psi_M$ has singular weight. This proves the theorem.

The expansion of $\Psi_M$ at the level $2p$ cusp shows
\begin{cp1}
The denominator function of the generalized Kac-Moody superalgebra obtained by twisting the fake monster superalgebra with an element of cycle shape $1^mp^m$ defines a holomorphic automorphic form of singular weight. 
\end{cp1}
The expansion of $\Psi_M$ at the other cusp implies
\begin{cp1}
There is a generalized Kac-Moody superalgebra with the following properties.
The root lattice is the dual $L'$ of the lattice $L=E_8^p(-1)\oplus II_{1,1}(2p)$. The Weyl group is the group generated by the reflections in the norm $-1/p$ vectors of $L'$ and the norm $-1$ vectors $\al$ in $L'$ with $2\al\in L$. The Weyl vector is $\rho=(0,0,1/2p)$. The real simple roots are the simple roots of $G$, i.e. the roots satisfying $(\rho,\al)=\al^2/2$. The imaginary simple roots are the positive multiples $n\rho$ of the Weyl vector with multiplicity $(-1)^n2m$ if $p$ divides $n$ and $(-1)^nm$ else. A root $\al$ is odd if and only if $p\al^2/2$ is integral and $\al$ has even order. The multiplicity of an even root $\al\in L'$ is $mult(\al) = c(p\al^2/2)$ if $2\al\notin L$ and $mult(\al) = c(\al^2/2) + c(p\al^2/2)$ if $2\al \in L$. For odd roots the same formula holds with opposite signs. The denominator identity is given by the expansion of $\Psi_M$ at the level $1$ cusp. 
\end{cp1}
Again these results can be extended to $p=1$.

\section{Hyperbolic reflection groups}

In this section we describe the reflection groups of the root lattices of the twisted fake monster algebras. In the bosonic case we get some information about these groups from the singular theta correspondence. In the super case we work out their fundamental domains using Vinberg's algorithm. 

\subsection{Lorentzian lattices}

Let $L$ be a Lorentzian lattice of dimension $n$. There are 2 cones of positive norm vectors in $L\otimes {\mathbb R}$. The vectors of norm $1$ in one of these cones form a copy of the $n-1$ dimensional hyperbolic space $H$. The automorphism group $Aut(L)$ of $L$ is the direct product of ${\mathbb Z}_2$ and the group $Aut(L)^+$ fixing the 2 cones of positive norm vectors. The reflection group $W$ of $L$ is the subgroup of $Aut(L)^+$ generated by reflections in the roots of $L$. W acts on $L\otimes {\mathbb R}$ and by restriction on $H$. The reflection hyperspaces divide $H$ into Weyl chambers. We choose one Weyl chamber $D$ and call it the fundamental Weyl chamber. Then $Aut(L)^+$ is the semidirect product $Aut(L)^+=Aut(D).W$ where $Aut(D)$ is the automorphism group of $D$. $W$ is called arithmetic if $Aut(D)$ is finite. The roots corresponding to the faces of $D$ form a set of simple roots of $L$. The reflections in these roots generate $W$. The angles between the simple roots and thus the defining relations of $W$ are usually described through the Dynkin diagram of $L$. 

Vinberg \cite{V} describes the following algorithm for finding a set of simple roots of $L$.
 
Choose a vector $w$ in $L$ with $w^2\geq 0$. The roots orthogonal to $w$ form a root system which is finite if $w^2>0$ and affine else. Choose a fundamental Weyl chamber $C$ for this root system. Then there is unique fundamental Weyl chamber $D$ of $W$ containing $w$ and contained in $C$, and its simple roots can be found inductively as follows. 

All the simple roots of $C$ are simple roots of $D$. Order the roots $\al$ which have positive inner product $(\al,w)$ with respect to the distance $(\al,w)/\sqrt{-\al^2}$ of their hyperplanes from $w$.

Take a root $\al$ as simple root for $D$ if and only if it has nonnegative inner product with all the simple roots we already have. It is sufficient to check this for the simple roots whose hyperplanes are strictly closer to $w$ than the hyperplane of $\al$. 

If at some point the roots we have already found by this algorithm span $L\otimes {\mathbb R}$ and contain at least one Dynkin diagram that is spherical of rank $n-1$ or affine of rank $n-2$, and every spherical diagram of rank $n-2$ in such a diagram is contained in a second such diagram, then these roots form a complete set of simple roots for $D$ (cf. for example Theorem 1.4 in \cite{B1}). 

\subsection{Reflection groups}

First we consider the root lattices of the fake monster algebras. They are similar to the lattice $II_{1,25}=\Lambda(-1) \oplus II_{1,1}$. Let $L=\Lambda^p(-1)\oplus II_{1,1}$. The roots of $L$ are the norm $-2$ vectors in $L$ and the norm $-2p$ vectors in $pL'$. The reflection group of $L$ is also the Weyl group of the fake monster algebra with root lattice $L$. Let $M=L\oplus II_{1,1}(p)$ and $f_{M+\delta}(\tau)$ be the components of the modular form $F$ with representation $\rho_M$ in section 3.2. We define 
\[ F_L(\tau)=\sum_{\gamma\in L'/L} f_{L+\gamma}(\tau)e^{\gamma} \]
with components
\[ f_{L+\gamma}(\tau) = \sum_{\delta \in M'/M \atop \delta|K=\gamma}
                        f_{M+\delta}(\tau)                           \]
where $z$ is a primitive norm $0$ vector in $II_{1,1,}(p)$ and 
$K=M\cap z^{\bot}$. It is easy to see that $F_L$ is a modular form of type $\rho_L$. The singular theta correspondence associates a modular form $\Phi_L$ to $F_L$ whose singularities are exactly at the reflection hyperplanes of $W$. Theorem \ref{st2} implies that the norm $0$ vector $\rho=(0,0,1)$ is Weyl vector for $L$ and the simple roots of $L$ are the roots $\al$ satisfying $(\rho,\al)=\al^2/2$. Furthermore the quotient $Aut(L)^+/W$ contains a free abelian subgroup of finite index. 

The lattice $\Lambda^p(-1)$ has no roots so that $Aut(L)^+/W$ is equal to the group of affine automorphisms of $\Lambda^p(-1)$ by Theorem 3.3 of \cite{B1}.

Now we consider the root lattices of the fake monster superalgebras. We will see that they are similar to the lattice $II_{1,9}=E_8(-1) \oplus II_{1,1}$. Here we apply Vinberg's algorithm rather than Theorem \ref{st2} because the latter would not give much information on the reflection groups. 

The lattice $E_8^3$ is the 4 dimensional lattice with elements $(m_1,m_2,m_3,m_4)$, where all $m_i$ are in $\mathbb Z$ or all $m_i$ are in $\mathbb Z+\frac{1}{2}$ and $\sum m_i$ is even, and norm $(m_1,m_2,m_3,m_4)^2=m_1^2+m_2^2+3m_3^2+3m_4^2$. The dual is the lattice with elements $(m_1,m_2,m_3/3,m_4/3)$ where the $m_i$ are as in $E_8^3$ and $(m_1,m_2,m_3/3,m_4/3)^2=m_1^2+m_2^2+m_3^2/3+m_4^2/3$. The norm 6 vectors of $E_8^3$ are all in $3{E_8^3}'$ so that the roots of $E_8^3$ are the norm 2 and norm 6 vectors. The root system of $E_8^3$ decomposes into 2 orthogonal $G_2$'s. 

We choose the vector $w=(0,1,-1)$ in $L=E_8^3(-1)\oplus II_{1,1}$ with norm $w^2=2$ and apply Vinberg's algorithm to determine the simple roots of $L$. We find the following complete set of simple roots
\renewcommand{\arraystretch}{1.3}
\[ \begin{array}{ll} 
\al_1=(v_1,0,0),  \qquad \quad &  
                     v_1=(\frac{1}{2},\frac{1}{2},-\frac{1}{2},-\frac{1}{2}),\\
\al_2=(v_2,0,0),  &  v_2=(0,0,1,1),                                          \\
\al_3=(v_3,0,0),  &  v_3=(\frac{1}{2},-\frac{1}{2},-\frac{1}{2},\frac{1}{2}),\\
\al_4=(v_4,0,0),  &  v_4=(0,0,1,-1),                                         \\
\al_5=(v_5,1,1),  &  v_5=(0,0,0,0),                                         \\
\al_6=(v_6,0,-1), &  v_6=(-1,1,0,0),                                        \\
\al_7=(v_7,0,-1), &  v_7=(-1,-1,0,0) .                                        
\end{array} \]
The roots $\al_1$ and $\al_4$ have norm $-6$ while the other simple roots have norm $-2$. $L$ also has a Weyl vector.
\begin{p1}
The Lorentzian lattice $L=E_8^3(-1)\oplus II_{1,1}$ has 7 simple roots which can be taken to be the roots $\al\in L$ with
\[ (\rho,\al)=\al^2/2 \]
where $\rho$ is the Weyl vector 
\[ \rho = (r,-6,7) \quad \mbox{with} \quad r=(5,0,1,0) \]
of norm $\rho^2=56$. The Dynkin diagram of $L$ is

\vspace{3ex}
\begin{figure}[h]
\begin{center}
\setlength{\unitlength}{0.00083333in}
\begingroup\makeatletter\ifx\SetFigFont\undefined%
\gdef\SetFigFont#1#2#3#4#5{%
  \reset@font\fontsize{#1}{#2pt}%
  \fontfamily{#3}\fontseries{#4}\fontshape{#5}%
  \selectfont}%
\fi\endgroup%
{\renewcommand{\dashlinestretch}{30}
\begin{picture}(3766,181)(0,-10)
\put(683,83){\ellipse{150}{150}}
\put(83,83){\ellipse{150}{150}}
\path(158,83)(608,83)
\path(153,128)(628,128)(623,128)
\path(148,33)(633,33)
\put(3683,83){\ellipse{150}{150}}
\put(3083,83){\ellipse{150}{150}}
\path(3158,83)(3608,83)
\path(3153,128)(3628,128)(3623,128)
\path(3148,33)(3633,33)
\put(1883,83){\ellipse{150}{150}}
\put(2483,83){\ellipse{150}{150}}
\put(1283,83){\ellipse{150}{150}}
\path(1358,83)(1808,83)
\path(1958,83)(2408,83)
\path(2558,83)(3008,83)
\path(758,83)(1208,83)
\end{picture}
}
\end{center}
\end{figure}
\noindent
The reflection group of $L$ is arithmetic and has index 2 in $Aut(L)^+$.  
\end{p1}

The root lattice of the other fake monster superalgebra corresponding to $p=3$ is the dual of $E_8^3(-1)\oplus II_{1,1}(6)$. We rescale the root lattice by $2p=6$ to obtain the even lattice $L=E_8^3(-2)\oplus II_{1,1}$. This lattice has level $6$ and roots of norm $-2,-4,-6$ and $-12$ in $L,L\cap 2L',L\cap 3L'$ and $6L'$. If we choose again $w=(0,1,-1)$ in Vinberg's algorithm we find 
\begin{pp1}
$L=E_8^3(-2)\oplus II_{1,1}$ has 8 simple roots which can be taken to be the roots $\al\in L$ with
\[ (\rho,\al)=\al^2/2 \]
where $\rho$ is the Weyl vector 
\[ \rho = (r,-6,7) \quad \mbox{with} \quad r=(5,0,1,0) \]
of norm $\rho^2=28$. The reflection group of $L$ has index 2 in $Aut(L)^+$ and the corresponding Dynkin diagram is
\end{pp1}
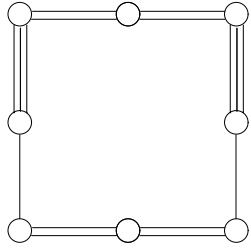
\begin{figure}[ht]
\begin{center}
\setlength{\unitlength}{0.00083333in}
\begingroup\makeatletter\ifx\SetFigFont\undefined%
\gdef\SetFigFont#1#2#3#4#5{%
  \reset@font\fontsize{#1}{#2pt}%
  \fontfamily{#3}\fontseries{#4}\fontshape{#5}%
  \selectfont}%
\fi\endgroup%
{\renewcommand{\dashlinestretch}{30}
\begin{picture}(1516,1531)(0,-10)
\put(83,1433){\ellipse{150}{150}}
\put(83,758){\ellipse{150}{150}}
\path(88,1348)(88,838)(88,828)
\path(88,1348)(88,838)(88,828)
\path(128,1368)(128,818)(128,823)
\path(48,1368)(48,818)(48,823)
\put(1433,1433){\ellipse{150}{150}}
\put(1433,758){\ellipse{150}{150}}
\path(1438,1348)(1438,838)(1438,828)
\path(1438,1348)(1438,838)(1438,828)
\path(1478,1368)(1478,818)(1478,823)
\path(1398,1368)(1398,818)(1398,823)
\put(83,83){\ellipse{150}{150}}
\put(758,83){\ellipse{150}{150}}
\put(1433,83){\ellipse{150}{150}}
\put(758,83){\ellipse{150}{150}}
\put(758,1433){\ellipse{150}{150}}
\put(758,1433){\ellipse{150}{150}}
\path(158,53)(688,53)(158,53)
\path(163,103)(688,103)(163,103)
\path(833,53)(1363,53)(833,53)
\path(838,103)(1363,103)(838,103)
\path(83,683)(83,158)
\path(1433,683)(1433,158)
\path(158,1403)(688,1403)(158,1403)
\path(163,1453)(688,1453)(163,1453)
\path(833,1403)(1363,1403)(833,1403)
\path(838,1453)(1363,1453)(838,1453)
\end{picture}
}
\end{center}
\end{figure}

$E^7_8$ is the 2 dimensional lattice with elements $(m_1,m_2)$ where either $m_1$ and $m_2$ are in $\mathbb Z$ and $m_1+m_2$ is even or $m_1$ and $m_2$ are in $\mathbb Z+\frac{1}{2}$ and $m_1+m_2$ is odd, and norm $(m_1,m_2)^2=m_1^2+7m_2^2$. The dual is the lattice $(m_1,m_2/7)$ with $m_1,m_2$ both in $\mathbb Z$ or both in ${\mathbb Z}+\frac{1}{2}$ and $m_1+m_2$ even. The norm in the dual lattice is $(m_1,m_2/7)^2=m_1^2+m_2^2/7$. The norm 14 vectors of $E^7_8$ are all in $7{E^7_8}'$ so that the roots of $E^7_8$ are the norm 2 and norm 14 vectors. The root system of $E^7_8$ is the sum of 2 orthogonal $A_1$'s where the roots in one $A_1$ are scaled to norm 2 and in the other $A_1$ to norm 14. While $E^3_8$ is generated by its short roots the roots of $E^7_8$ generate a sublattice of index 2.  

Again we use Vinberg's algorithm to determine the simple roots of $L=E_8^7(-1)\oplus II_{1,1}$. We choose $w$ as above. Then we find the following simple roots 
\renewcommand{\arraystretch}{1.3}
\[ \begin{array}{ll} 
\al_1=(v_1,0,0),  \qquad \quad  
                  &  v_1=(\frac{1}{2},\frac{1}{2}),\\
\al_2=(v_2,0,0),  &  v_2=(\frac{7}{2},-\frac{1}{2}),\\
\al_3=(v_5,1,1),  &  v_3=(0,0),                 \\ 
\al_4=(v_4,0,-1), &  v_4=(-\frac{1}{2},-\frac{1}{2}),\\
\al_5=(v_5,1,-1), &  v_5=(-2,0), \\ 
\al_6=(v_6,0,-7), &  v_6=(-\frac{7}{2},\frac{1}{2}).
\end{array} \]
The roots $\al_2$ and $\al_6$ both have norm $-14$ while the other simple roots have norm $-2$. We get   
\begin{pp1}
The Lorentzian lattice $L=E_8^7(-1)\oplus II_{1,1}$ has 6 simple roots and a Weyl vector of norm 8. The above roots correspond to the Weyl vector $\rho = (r,-2,3)$ with $r=(2,0)$. The reflection group of $L$ has index 2 in $Aut(L)^+$. 
The Dynkin diagram of $L$ is
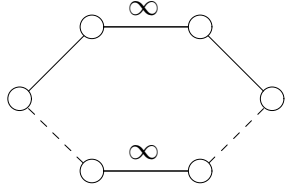
\vspace{1ex}
\begin{figure}[ht]
\begin{center}
\setlength{\unitlength}{0.00083333in}
\begingroup\makeatletter\ifx\SetFigFont\undefined%
\gdef\SetFigFont#1#2#3#4#5{%
  \reset@font\fontsize{#1}{#2pt}%
  \fontfamily{#3}\fontseries{#4}\fontshape{#5}%
  \selectfont}%
\fi\endgroup%
{\renewcommand{\dashlinestretch}{30}
\begin{picture}(1741,1217)(0,-10)
\put(1658,533){\ellipse{150}{150}}
\put(1208,83){\ellipse{150}{150}}
\put(1208,983){\ellipse{150}{150}}
\path(1133,83)(608,83)
\path(1258,928)(1603,583)
\path(1133,983)(608,983)
\dashline{60.000}(1603,473)(1263,133)
\put(83,533){\ellipse{150}{150}}
\put(533,83){\ellipse{150}{150}}
\put(533,983){\ellipse{150}{150}}
\path(608,83)(1133,83)
\path(483,928)(138,583)
\path(608,983)(1133,983)
\dashline{60.000}(138,473)(478,133)
\put(758,1058){\makebox(0,0)[lb]{\smash{{{\SetFigFont{12}{14.4}{\rmdefault}{\mddefault}{\updefault}$\infty$}}}}}
\put(758,158){\makebox(0,0)[lb]{\smash{{{\SetFigFont{12}{14.4}{\rmdefault}{\mddefault}{\updefault}$\infty$}}}}}
\end{picture}
}
\end{center}
\end{figure}
\end{pp1}

When we rescale the root lattice of the other fake monster superalgebra corresponding to $p=7$ by $14$ we obtain the lattice $L=E_8^7(-2)\oplus II_{1,1}$ of level $14$. This lattice has roots of norm $-2,-4,-14$ and $-28$ in $L, L\cap 2L', L\cap 7L'$ and $14L'$. We have
\begin{pp1}
The Lorentzian lattice $L=E_8^7(-2)\oplus II_{1,1}$ has 8 simple roots and a Weyl vector of norm 4. The Dynkin diagram of $L$ is
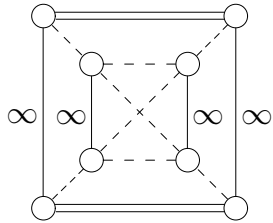
\vspace{1ex}
\begin{figure}[h]
\begin{center}
\setlength{\unitlength}{0.00083333in}
\begingroup\makeatletter\ifx\SetFigFont\undefined%
\gdef\SetFigFont#1#2#3#4#5{%
  \reset@font\fontsize{#1}{#2pt}%
  \fontfamily{#3}\fontseries{#4}\fontshape{#5}%
  \selectfont}%
\fi\endgroup%
{\renewcommand{\dashlinestretch}{30}
\begin{picture}(1814,1381)(0,-10)
\put(1430,1283){\ellipse{150}{150}}
\put(230,1283){\ellipse{150}{150}}
\path(305,1308)(1360,1308)
\path(310,1263)(1360,1263)
\put(1430,83){\ellipse{150}{150}}
\put(230,83){\ellipse{150}{150}}
\path(305,108)(1360,108)
\path(310,63)(1360,63)
\put(530,383){\ellipse{150}{150}}
\put(1130,383){\ellipse{150}{150}}
\put(530,983){\ellipse{150}{150}}
\put(1130,983){\ellipse{150}{150}}
\path(1430,1208)(1430,158)
\path(230,1208)(230,158)
\path(530,908)(530,458)
\path(1130,908)(1130,458)
\dashline{60.000}(605,983)(1055,983)
\dashline{60.000}(605,383)(1055,383)
\dashline{60.000}(285,1223)(470,1038)
\dashline{60.000}(1195,323)(1380,138)
\dashline{60.000}(590,923)(1075,438)
\dashline{60.000}(1080,933)(585,438)
\dashline{60.000}(1380,1228)(1185,1033)
\dashline{60.000}(480,328)(285,138)
\put(305,608){\makebox(0,0)[lb]{\smash{{{\SetFigFont{12}{14.4}{\rmdefault}{\mddefault}{\updefault}$\infty$}}}}}
\put(0,613){\makebox(0,0)[lb]{\smash{{{\SetFigFont{12}{14.4}{\rmdefault}{\mddefault}{\updefault}$\infty$}}}}}
\put(1430,608){\makebox(0,0)[lb]{\smash{{{\SetFigFont{12}{14.4}{\rmdefault}{\mddefault}{\updefault}$\,\infty$}}}}}
\put(1125,608){\makebox(0,0)[lb]{\smash{{{\SetFigFont{12}{14.4}{\rmdefault}{\mddefault}{\updefault}$\,\infty$}}}}}
\end{picture}
}
\end{center}
\end{figure}

\noindent 
The reflection group of $L$ is arithmetic and has index 2 in $Aut(L)^+$.  
\end{pp1} 
We remark that $L$ has 2 simple roots of each possible root length. The Weyl vector corresponding to the choice $w=(0,0,1)$ in Vinberg's algorithm is $\rho = (r,-2,3)$ with $r=(2,0)$. 

\section*{Acknowledgments}
I thank R. E. Borcherds for stimulating discussions.

\end{document}